\documentstyle[12pt,amsfonts,amscd]{amsart}
\newtheorem{Theorem}{Theorem}[section]
\newtheorem{Definition}[Theorem]{Definition}
\newtheorem{Proposition}[Theorem]{Proposition}

\newtheorem{Lemma}[Theorem]{Lemma}
\newtheorem{Corollary}[Theorem]{Corollary}
\theoremstyle{remark}
\newtheorem{Example}[Theorem]{Example}

\def\il{\int\limits_}
\def\CC{{\Bbb C}}
\def\RR{{\Bbb R}}

\def\cara{Carath\'eodory }
\def\epf{\hskip.2in\vrule width.4pt height6pt depth0pt\vrule
width5.2pt height6pt depth-5.6pt\hskip-5.2pt\vrule width5.2pt
height.4pt depth0pt\vrule width.4pt height6pt depth0pt\ }

\def\LL{\Lambda}

\def\eps{\varepsilon}

\def\ovr{\overline}

\def\pd{\partial}

\def\Dl{\Delta}
\def\dl{\delta}

\def\bd{\partial}

\def\sbs{\subset}

\def\Sbs{\subset\subset}

\def\wtl{\widetilde}

\def\Aut{\operatorname{Aut}}

\def\re{{\mathbf {Re\,}}}
\def\im{\mathbf {Im\,}}
\def\be{\begin{enumerate}}
\def\ee{\end{enumerate}}
\def\bT{\begin{Theorem}}
\def\eT{\end{Theorem}}
\def\bP{\begin{Proposition}}
\def\eP{\end{Proposition}}
\def\bD{\begin{Definition}}
\def\eD{\end{Definition}}
\def\bE{\begin{Example}}
\def\eE{\end{Example}}
\def\bL{\begin{Lemma}}
\def\eL{\end{Lemma}}
\def\bC{\begin{Corollary}}
\def\eC{\end{Corollary}}

\def\H{{\mathcal H}}

\begin{document}
\title{ Upper semicontinuity of the dimensions of automorphism groups
of domains in $\CC^n$ }
\author{Buma L. Fridman, Daowei Ma and Evgeny A. Poletsky}
\begin{abstract} Let $\H^n$ be the metric space of all bounded domains
in ${\Bbb C}^n$ with the metric equal to the Hausdorff distance
between boundaries of domains. We prove that the dimension of the
group of automorphisms of domains is an upper semicontinuous function
on $\H^n$. We also provide theorems and examples regarding the
change in topological structure of these groups under small
perturbation of a domain in $\H^n$.\end{abstract}
\thanks{E.A. Poletsky was partially supported by NSF Grant
DMS-9804755.}
\address{ buma.fridman@@wichita.edu, Department of Mathematics,
Wichita State University, Wichita, KS 67260-0033, USA}
\address{ dma@@math.twsu.edu, Department of Mathematics,
Wichita State University, Wichita, KS 67260-0033, USA}
\address{ eapolets@@syr.edu, Department of Mathematics,  215
Carnegie Hall, Syracuse University,  Syracuse, NY 13244, USA}
\maketitle \setcounter{section}{-1}
\section{Introduction}
The automorphism group $\Aut(D)$ (the group of biholomorphic self-maps
of $D$) of a bounded domain $D$ in $\CC^n$ is, in general, difficult to
describe and little is known about it. However, it is known (see
\cite{SZ,BD}) that any compact Lie group can be realized as the group
of automorphisms of a smooth strictly pseudoconvex domain, and (see
\cite{ShT}) that any linear Lie group can be realized as the group of
automorphisms of a bounded domain. So, if we consider the group
$\Aut(D)$ as a function of $D$, the set of values is quite large.
\par If one considers this function on the metric space $\H^n$ of all
bounded domains in ${\Bbb C}^n$ with the metric equal to the
Hausdorff distance between boundaries of domains, one can expect that
small perturbation of the boundary may only ``decrease'' the group,
i.e., the function $\Aut(D)$ is ``upper semicontinuous''. Indeed, in
[GK], [Ma] and [FP] the authors, using topologies different from
$\H^n$, proved the upper semicontinuity of the function $\Aut(D)$ in
the sense that $\Aut(\wtl D)$ is isomorphic to a subgroup of
$\Aut(D)$ when $\wtl D$ is ``close'' to $D$. But, in general, this
idea is not true according to the following theorem (\cite{FP}).
\bT\label{T:fp}
Let $M$ be a domain in $\CC^n$. Then there exists an increasing
sequence of bounded domains $M_k\subset M_{k+1}\subset\subset M$ such
that $M=\cup M_k$ and $\Aut(M_k)$ contains a subgroup isomorphic to
$Z_k$. \eT
\par This shows that domains in $\CC^n$ with an automorphism group
containing $Z_k$  are everywhere dense in $\H^n$,  and it is well
known that domains without non-trivial automorphisms are dense in
$\H^n$. So arbitrarily small perturbation of a domain in $\H^n$ may
create a domain with a larger automorphism group. But, for all known examples,
this group is
discrete, so it is of dimension zero. The natural question arises:
can small perturbation in $\H^n$ create domains with larger
dimensions of automorphism groups?
\par In this paper we answer this question in the negative. Namely,
we prove the following
\bT\label{T:mt} The function $\dim\Aut(D)$ is upper semicontinuous on
$\H^n$.
\eT
An immediate consequence is the following
\bC\label{C:gc} For
each $k>0$ the set of all domains in $\H^n$ whose groups of
automorphisms have dimensions greater than or equal to $k$ is closed
and, therefore, nowhere dense.\eC
Thus a domain cannot be approximated by domains whose automorphism groups
have strictly larger dimensions.
\par To prove Theorem \ref{T:mt} we consider a sequence of domains
$D_j$ converging in $\H^n$ to a domain $D$. The identity components
$\Aut_0(D_j)$ of $\Aut(D_j)$ have the same dimensions as $\Aut(D_j)$.
Also the dimensions of the Lie algebras of holomorphic vector fields
generated by all one-parameter groups in $\Aut_0(D_j)$ coincide with
$\dim\Aut_0(D_j)$. Lemma \ref{L:ne} states that the uniform norm of
such fields on a compact set is bounded by its norm on an arbitrarily
selected ball times a constant that, basically, depends on the size
of the ball and the distances from the ball and the compact set to the
boundary of a domain. This allows us to normalize bases in Lie
algebras of $\Aut_0(D_j)$ and apply Theorem~\ref{T:lt}, which asserts
the existence of non-trivial limits of those vector fields. The
limits belong to the Lie algebra of $\Aut_0(D)$ and this gives us the
proof.
\par It is reasonable to ask whether $\Aut_0(D_j)$ are always isomorphic to
a subgroup of $\Aut_0(D)$ when $j$ is large. An example in Section
\ref{S:st} shows that the answer is negative.

If $K_j$ is a maximal compact subgroup of $\Aut_0(D_j)$, then
$\Aut_0(D_j)$ is diffeomorphic to $K_j\times{\Bbb R}^{k_j}$ (see [MZ,
p.~188]). The groups $K_j$ may decrease or even  disappear in the limit
(see Example \ref{E:1}), while non-compact parts  never vanish (see
Theorem \ref {T:nc}).

\section{Some basic facts}

\par Let $D$ be a bounded domain in $\CC^n$. If the Lie group $\Aut(D)$
has positive dimension, then it has one-parameter subgroups $g(\cdot,t)$,
$-\infty<t<\infty$, i.e., $g(z,t+s)=g(g(z,t),s)$. Such subgroups
generate vector fields $$X(z)=\frac{\bd g}{\bd t}(z,0)$$ that are
holomorphic. Also, if $X$ is a holomorphic vector field on $D$ that
is ${\Bbb R}$-complete, i.e., the initial value problem
$$\frac{\partial g}{\partial t}(z,t)=X(g(z,t)), \;\;g(z,0)=z$$
has a solution on $D\times{\Bbb R}$, then
$g(z,t)$ is a one-parameter group.
\par The vector field $X$ has the following group property:
\begin{equation}\label{e:gp}
X(g(z,t))=\frac{\bd g}{\bd z}(z,t)X(z).
\end{equation}
\par For every two points $z$ and $w$ in $D$ among all holomorphic
mappings of $D$ into the unit disk $\Dl$ we choose  holomorphic
functions $f$ such that $f(w)=0$ and $f(z)$ is real and the maximal
possible. Such functions $f$ exist and are called \cara extremal
functions for $z$ and $w$ on $D$. The quantity
\begin{equation}\label{e:cm}
\rho(0,f(z))=\frac12\ln\frac{1+f(z)}{1-f(z)}
\end{equation}
is called the \cara distance $c_D(z,w)$ on $D$. (Note that the
formula for $\rho(0,a)$ gives the Poincar\'e distance between $0$
and $a$ in the unit disc.) When $D$ is bounded this distance is
non-degenerate and invariant, i.e., $c_D(g(z),g(w))=c_D(z,w)$ for
every $g\in \Aut(D)$ (see \cite[Ch.~5, \S 18]{Sh}).
\par For a point $w\in D$ and a vector $Y$ in $\CC^n$, among all
holomorphic mappings of $D$ into the unit disk $\Dl$  we choose
holomorphic functions $f$ such that $f(w)=0$ and $(f'(w),Y)$ is real
and the maximal possible. (Here $(Z,Y)=\sum_{j=1}^nz_jy_j$.) These
functions are \cara extremal functions for $Y$ at $w$ in $D$. The
Carath\'eodory length of $Y$ at $w$ is defined to be $C_D(w,Y)=(f'(w),Y)$.
It follows from \cite[Ch.~5, \S 18]{Sh} that if $w(t)$ is a smooth curve
in $D$ with $w(0)=w$ and $Y=w'(0)$, then
\begin{equation}\label{e:cn}
c_D(w,w(t))=C_D(w,Y)t+o(t).
\end{equation}
Let $B(w,r)$ be the ball of radius $r$ centered at $w$ and let $|Y|$
be the Euclidean norm of $Y$. If $B(w,r)\sbs D\sbs B(w,R)$, then
\begin{equation}\label{e:cne}\frac{|Y|}{R}\le C_D(w,Y)\le\frac{|Y|}{r}.
\end{equation}
\section{Proof of Theorem 0.2}\label{S:pmt}
\bL\label{L:id}
Let $D$ be a domain in $\CC^n$ and let $d(z,w)$ be an invariant
metric on $D$ satisfying the triangle inequality. If $g(z,t)$ is a
group action on $D$, then for any $w, z\in D$
$$|d(g(w,t),z)-d(w,z)|\le d(z, g(z,t)).$$
\eL
\begin{pf} Apply the identity $d(w,z)=d(g(w,t), g(z,t))$
and the triangle inequality.
\end{pf}
\bL\label{L:de}
Let $w\in B(w, r)\Sbs D\Sbs B(w, R)\Sbs \CC^n$. Then for any $Y\in
\CC^n$, $|Y|=1$, $$\re(\nabla f_s(w),Y)>\frac1{4R},$$ where $f_s(z)$
is a \cara extremal function for $w$ and $w+sY$ in $D$, and $s$ is a
real number such that
$$0<s\le\eps=\frac{r^2}{16R}.$$
\eL
\begin{pf} Let us fix $Y$ and introduce
$D_Y=\{\xi\in \CC: w+\xi Y\in D\}$.
Clearly, $\Dl(0, r)\Sbs D_Y \Sbs \Dl(0,R)$, where $\Dl(0,s)$ is the
disk of radius $s$ centered at 0.
\par Let $g_Y(z)$ be a \cara extremal function for $Y$ at $w$ in $D$.
For $\xi\in D_Y$ we introduce the functions $u(\xi)=\re
F(\xi)$, where $F(\xi)=f_s(w+\xi Y)$, and $v(\xi)=\re G(\xi)$, where
$G(\xi)=g_Y(w+\xi Y)$. All these functions are well-defined on $D_Y$
and $F(0)=G(0)=0$, $G'(0)= C_D(w,Y)$ and $u(s)=F(s)\geq v(s)$. Let us
prove that $v(t)\geq t/(2R)$ when $t\in[0,\eps]$. Since
$$|v(t)-v'(0)t |\leq\frac12\sup_{0\le x\leq\eps} |v''(x)|\cdot t^2,$$
$\eps\le r^2/(16R)<r/2$ and by Cauchy estimate
$|v''(x)|\le 2/(r-\eps)^2$ when $x<\eps$, we see that
$$v(t)\ge v'(0)t-\frac{1}{(r-\eps)^2}t^2\ge
v'(0)t-\frac4{r^2}t^2.$$
Since $t\le r^2/(16R)$ and by (\ref{e:cne})
$$v'(0)=G'(0)=C_D(w,Y)\ge\frac1R,$$
$$v(t)\ge \frac tR-\frac t{4R}\ge
\frac t{2R}.$$
for $0\le t\le\eps$.
In particular, $$v(s)=\re G(s)\ge\frac{s}{2R}.$$
\par Applying to the function $u(t)$ the same analysis as above we
obtain
$$u(s)\le u'(0)s+\frac{4}{r^2}s^2.$$
Hence
$$\frac s{2R}\le v(s)\le u(s)\le u'(0)s+\frac4{r^2}s^2\le
u'(0)s+\frac{s}{4R}.$$
Thus
$$\re(\nabla f_s(w),Y)=\re F'(0)=u'(0)\ge\frac1{4R}.$$
\end{pf}
\bL\label{L:one}
Let $B(0,r+a)\Sbs D\Sbs B(0,R)$, $r,a>0$. Then there exists a positive
$\delta=\delta(a,r,R)<a$ such that
$$\|X\|_{B(0,r+\delta)}\le \frac{32R}{a}\|X\|_{B(0,r)}$$
for every holomorphic vector field $X$ generated by
a one-parameter group action $g(z,t)$ on $D$.\eL
\begin{pf}
Let $w$ belong to $B(0,r+a/2)$. Since
$$w\in B(w, a/2)\Sbs D\Sbs B(w, 2R),$$
by Lemma \ref{L:de} there is an $\eps=\eps(a,R)>0$ such that for every
$w\in B(0, r+a/2)$, every $Y\in\CC^n$,
$|Y|=1$, and every $s\in(0,\eps]$
\begin{equation}\label{e:1}
\re(\nabla f(w), Y)\ge \frac1{8R},
\end{equation}
where $f$ is a \cara extremal function for $w$ and $w+sY$.
\par Let us take a positive number $\dl<a/2$ so small that for every
$w\in B(0,r+\dl)$ and every unit vector $V$ there is a unit vector
$Y$ such that $w+sY\in B(0,r)$ for some real $s$ with $|s|<\eps$ and
$$|V-Y|<b=\frac{a}{32R}.$$ Clearly, the choice of this $\dl$ depends
only on $a$, $r$ and $R$.
\par The lemma needs a proof only for non-trivial group actions when
$X\not\equiv0$.
Let $w\in\bd B(0, r+\dl)$, $X(w)\ne0$ and let $V=X(w)/|X(w)|$. We
choose a vector $Y$ and a real $s$ satisfying the above conditions.
\par Let $f$ be a \cara extremal function for $w$ and $z=w+s Y$. Since
$B(w,a/2)\Sbs D$, by Schwarz inequality,
$$|\re(\nabla f(w), Y-V)|\le b|\nabla f(w)|\le \frac{2b}a.$$
Hence by (\ref{e:1}),
$$\re(\nabla f(w), V)\ge \re(\nabla f(w), Y)-\frac{2b}a\ge
\frac1{8R}-\frac{1}{16R}=\frac1{16R}.$$
\par Let $\zeta(t)=f(g(w,t))$ and $p=f(z)$. We introduce
$$m(t)=\left|\frac{\zeta(t)-p}{1-p\zeta(t)}\right|.$$
If $\rho(\zeta,\xi)$ is the Poincare metric on $U$, then
$\rho(0,p)=c_D(z,w)$ and
$$\rho(\zeta(t),p)=\frac12\ln\frac{1+m(t)}{1-m(t)}.$$
A straightforward calculation shows that
$$\frac{dm^2}{dt}(0)=-2p(1-p^2)\re (\nabla f(w),X(w)),$$
and, by using this calculation, we obtain
$$\frac{d}{dt}\rho(\zeta(0),p)=-\re(\nabla f(w),X(w)).$$
Hence
$$\rho(\zeta(-t),p)\ge\rho(0,p)+t\re (\nabla f(w),X(w))
\ge c_D(z,w)+\frac{t}{16R}|X(w)|$$
for small positive $t$. Since the \cara
metric decreases under the holomorphic mapping $f$,
\begin{equation}\label{e:2}
c_D(z,g(w,-t))\ge\rho(\zeta(-t),p)\ge c_D(z,w)+\frac{t}{16R}|X(w)|.
\end{equation}
By (\ref{e:2}) and Lemma~\ref{L:id},
$$\frac{t}{16R}|X(w)|\le c_D(z,g(w,-t))-c_D(z,w)
\le c_D(z,g(z,-t)).$$
By (\ref{e:cn}),
$c_D(z,g(z,-t))=C_D(z,X(z))t+o(t)$.
Note that $B(z,a)\sbs D$ and, therefore, $C_D(z,X(z))\le1/a$. Hence
$$c_D(z,g(z,-t))\le 2C_D(z,X(z))t\le\frac2a|X(z)|t$$ for small positive
$t$. Thus $$|X(w)|\le \frac{32R}{a}|X(z)|$$  and
$$\|X\|_{B(0,r+\dl)}\le \frac{32R}{a}\|X\|_{B(0,r)}.$$\end{pf}

\bL\label{L:ne} Let $R>2r>2s>0$. Let $K$ be a connected compact set
containing 0 in $\CC^n$. Let $D$ be a domain in $\CC^n$ such that
$B(0,2r)\sbs D\sbs B(0,R)$ and such that the $3s$-neighborhood of $K$
is contained in $D$. Then there exists a positive constant $C=C(K, R,
s)$ such that $\|X\|_K\le C\|X\|_{B(0,r)}$ for each holomorphic
vector field $X$ generated by a one-parameter group action $g(z,t)$
on $D$. \eL
\begin{pf}
Let $X$ be such a vector field on $D$. By the previous lemma there exist
positive numbers  $\delta=\delta(s,R)<s$ and $c=c(s,R)$ such that
$$\|X\|_{B(z, s+\delta)}\le c\|X\|_{B(z,s)}$$
whenever $z\in D$ is at least $3s$ away from $\pd D$.
There is a positive integer $N=N(K, \delta)$
such that for each $z\in K$ there is a set of $N$ points
$\{z_1,\dots, z_N\} \sbs K$
with $z_1=0$, $z_N=z$, and $|z_{k+1}-z_k|<\delta$ for $k=1,\dots, N-1$.
Since $B(z_{k+1},s)\sbs B(z_k, s+\delta)$, we see that
$$\|X\|_{B(z_{k+1}, s)}\le c\|X\|_{B(z_{k}, s)}$$
for $k=1, \dots, N-1$. Thus,
$$\|X\|_{B(z, s)}\le c^{N-1}\|X\|_{B(0, s)}.$$
In particular, $|X(z)|\le c^{N-1}\|X\|_{B(0,s)}\le c^{N-1}\|X\|_{B(0,r)}$.
Therefore, $\|X\|_K\le c^{N-1}\|X\|_{B(0,r)}$.
\end{pf}

\bT\label{T:lt} Suppose a sequence of domains $D_j$
converge in $\H^n$ to a domain $D$ and a ball $B(p,r+a)$, $r,a>0$,
belongs to all $D_j$. Also suppose that $g_j(z,t)$ are non-trivial
one-parameter group actions on $D_j$ generating the holomorphic
vector fields $X_j$. If $\|X_j\|_B=1$, $B=B(p,r)$, then there is a
subsequence of the group actions $g_{j_k}(z,t)$ that converges to a
non-trivial group action $g(z,t)$ on $D$ uniformly on compacta in
$D\times \RR$ and
$$\lim_{k\to\infty}X_{j_k}(w)=X(w)$$ uniformly on compacta in $D$, 
where $X$ is the holomorphic vector field generated by $g$.
\eT
\begin{pf} Let $K\Sbs D$. Choose $\delta>0$ so that the
$3\delta$-neighborhood of $K$ is contained in $D$ and in each $D_j$.
Let $\tilde K$ and $\hat K$ denote the $2\delta$-neighborhood and the
$\delta$-neighborhood of $K$ respectively. By Lemma \ref{L:ne} there
exists $A>0$ such that $\|X_j\|_{\wtl K}\le A$.
\par Let $\tau=\dl/(2A)$. Define the mapping
$h_j:\,\hat K\times(-\tau,\tau)\to D_j$ as the solution of the initial
value problem
$$\frac{\bd}{\bd t}h_j(z,t)=iX_j(h_j(z,t)),\;\;\;h_j(z,0)=z.$$
Since $\tau|X_j|<\dl$ in $\tilde K$, it follows from the ODE's theory
that the mapping $h_j$ is well-defined.
\par For $M=\{\zeta\in{\Bbb C}:\,|\im\zeta|<\tau\}$ we define
$G_j:\hat K\times M\to D_j$ by $G_j(z,t+is)=g_j(h_j(z,s),t)$. Since
$X_j$ is holomorphic, the mapping $G_j$ is holomorphic in $z$. We now
prove that it is holomorphic in $\zeta=t+is$. It is clear that
\begin{equation}\label{e:cr1}
\frac{\bd G_j}{\bd t}(z,t+is)=X_j(G_j(z,t+is)).
\end{equation}
It follows immediately from the fact that the Poisson brackets
$[X_j,iX_j]\equiv0$, that
\begin{equation}\label{e:cr2}
\frac{\bd G_j}{\bd s}(z,t+is)=iX_j(G_j((z,t+is)).
\end{equation}
This fact also can be proved by a straightforward reasoning:
\begin{equation}\begin{align}
\frac{\bd G_j}{\bd s}(z,t+is)&
=\frac{\bd g_j}{\bd z}(h_j(z,s),t)\cdot iX(h_j(z,s))\notag\\
&=iX_j(g_j(h_j(z,s),t))=iX_j(G_j(z,t+is));\notag
\end{align}\end{equation}
the middle equality is by the infinitesimal group property
(\ref{e:gp}). The equations (\ref{e:cr1}) and (\ref{e:cr2}) are the
Cauchy-Riemann equations for $G_j$ in $\zeta$. So $G_j$ is
holomorphic.
\par Passing to a subsequence, if necessary,  we may assume that the
mappings $G_j$ converge to a mapping $G$  uniformly on compacta in
$\hat K\times M$. Consequently, the mappings $g_j(z,t)$ converge to
$g(z,t)$ uniformly on compacta in $\hat K\times{\Bbb R}$, and the
vector fields $X_j$ converge to
$$X(z)=\frac{\bd g}{\bd t}(z,0)$$
uniformly on compacta in $\hat K$.
\par It follows that some subsequence of the sequence $\{g_j(z,t)\}$ converges to a mapping $g(z,t)=G(z,t)$ uniformly on  compacta in $D\times{\Bbb R}$.
Thus, $g(z,t)$ is a group action. Since $\|X\|_B=1$, this group
action is non-trivial.
\end{pf}

\noindent {\it Proof of Theorem \ref{T:mt}.} Let $D_j$ be a sequence
of domains converging in $\H^n$ to a domain $D$. Let  us choose
a ball $B(p,r+a)$, $r,a>0$, belonging to all $D_j$ for sufficiently
large $j$ and take $\dl>0$ from Lemma \ref{L:one}. Let $B=B(p,r)$ and
$\hat B=B(p, r+\delta)$. We may assume that the dimensions of all
groups $G_j=\Aut_0(D_j)$ are the same and equal to $k$. Since the Lie
algebra $A_j$ of all holomorphic vector fields on $D_j$ generated by
one-parameter subgroups in $G_j$ has the same dimension as $G_j$, we
can choose $X^m_j\in A_j$, $1\le m\le k$, such that
$$\il{\hat B}(X^m_j,\ovr X^l_j)\,dV=\dl_{ml},$$
where $\delta_{ml}$ is Kronecker's delta.
\par Clearly, $\|X_j^m\|_{\hat B}\ge \text{Vol}(\hat B)^{-1}$. On the other
hand, by Cauchy estimates and Lemma \ref{L:ne}, for some constants we have
$$1\ge C_1\|X_j^m\|_{B}\ge C_2\|X_j^m\|_{\hat B}.$$
\par Let $g^m_j$ be the one-parameter groups generated by $X^m_j$. By
Theorem \ref{T:lt} one can choose a subsequence  $\{j_k\}$ such that
$g^m_{j_k}$ converge, uniformly on compacta in $D\times \RR$, to a
one-parameter group $g^m(z,t)$ on $D$, and $X^m_{j_k}$ converge to a
vector field $X^m$ uniformly on compacta in $D$. Since
$$\il {\hat B}(X^m,\ovr X^l)\,dV=\dl_{ml},$$
the dimension of $\Aut_0(D)$ is at least $k$.

\section{Structural theorems}\label{S:st}
By Iwasawa's theorem (see [MZ, p.~188]) the group $\Aut_0(D)$ is
diffeomorphic to $K\times{\Bbb R}^{k}$, where $K$ is a maximal compact
subgroup and $k$ is the {\it characteristic number} of $\Aut(D)$. It
is interesting to find out what happens with $K$ and ${\Bbb R}^k$
under small perturbations of domains. Let us look at maximal compact
subgroups first. The argument of Corollary~4.1 in \cite{FP} provides
the following theorem.
\bT
Let $D$ be a bounded domain in ${\Bbb C}^n$, let $z_0$ be a point in
$D$, and let $W$ be a compact set in $D$. If $\wtl D$ is sufficiently close to $D$ in
$\H^n$ and for some maximal compact subgroup $\wtl K$ in $\Aut_0(\wtl
D)$ the orbit $\wtl K(z_0)\sbs W$, then $\wtl K$ is isomorphic to a
subgroup of $\Aut_0(D)$. \eT
\par Next example shows that without the condition in the above theorem
of orbits being contained in a fixed compact set, it is possible that
$\Aut(D)$ does not contain a compact subgroup while close domains
have $\Aut_0(\wtl D)$ isomorphic to $S^1$. Let $\Dl$ denote the unit
disc in $\CC$.

\bE\label{E:1} There is a sequence $\{D_j\}$ of bounded pseudoconvex
domains in $\CC^2$ converging to a domain $D$ such that
$\Aut(D_j)\cong  S^1$ for each $j$, and $\Aut(D)\cong\RR$.
\eE

\noindent {\it Construction.} Let $Q_j=\{z\in \Dl: |z- 2^{-1}+2^{-j}|>1/2\}$,
$Q=\{z\in \Dl: |z- 2^{-1}|>1/2\}$, $D_j=\{(z,w): z\in Q_j, w\in \Dl,
w\ne z\}$,  $D=\{(z,w): z\in Q, w\in \Dl, w\ne z\}$.

1. One can see that $D_j\to D$.

2. The domains $D_j$ and $D$ are bounded and pseudoconvex.

3. We now prove that $\Aut(D)\cong\RR$. Let $F\in \Aut(D)$. On each
fiber $(z,\cdot)$, $F$ is bounded and has an isolated singularity, so
$F$ extends to be an automorphism of $Q\times \Dl$. Thus, $F$ has the
form $F(z,w)=(f(z), g(w))$, or $F(z,w)=(g(w), f(z))$. For both cases,
one has, by the definition of $D$, that
\begin{equation}\label{e:fg}
f(z)=g(z),\;\;\;  z\in Q.\end{equation}
The second case is impossible, since implies that $f(Q)= \Dl$, $g(\Dl)=Q$, and
$f(Q)=g(Q)$, which leads to a contradiction that $\Dl$ coincides with a
subset of $Q$. Therefore, $F$ has the form $F(z,w)=(f(z), g(w))$,
where $f\in\Aut(Q)$, $g\in\Aut(\Dl)$. By (\ref{e:fg}), $f=g|_Q$. Let
$\phi(w)=- i(w+1)/(w-1)$. Then $\phi$ is a biholomorphic map from
$\Dl$ to the upper half-plane $\Pi=\{\zeta\in \CC: \im\zeta>0\}$, and
$\phi(Q)=\LL\equiv\{\zeta\in \CC: 0<\im\zeta<1\}$. Now $\phi\circ
g\circ \phi^{-1}$ is an automorphism of $\Pi$, and its restriction to
$\LL$ is an automorphism of $\LL$. Thus $\phi\circ g\circ
\phi^{-1}(\zeta)=\zeta+t$ for some $t\in \RR$. It follows that
$\Aut(D)=\{F_t: t\in \RR\}\cong \RR$, where
$F_t(z,w)=(g_t(z),g_t(w))$, and
$$g_t(w)=\phi^{-1}(\phi(w)+t)={2w+i(w-1)t\over 2+i(w-1)t}.$$
\par 4. In a way very similar to the above argument, one can prove that
$\Aut(D_j)\cong  S^1$ for each $j$. \epf
\par By Theorem \ref{T:mt} the creation of compact subgroups with
larger dimensions by small perturbations must be compensated by an
elimination of some non-compact subgroups so that the total dimension
will not go up. It seems to us that the other way around is
impossible: characteristic numbers are upper semicontinuous on
$\H^n$. While we cannot prove this statement, the following theorem
certifies that non-compact parts cannot be created from nothing.

\bT\label{T:nc} Let $D\sbs \CC^n$ be a bounded domain such that
$\Aut_0(D)$ is compact. Then for all $\wtl D$ sufficiently close in
$\H^n$ to $D$ the group $\Aut_0(\wtl D)$ is also compact. \eT
\begin{pf} If the statement is not true, then there is a sequence
$\{D_j\}$ of domains converging to $D$ such that for each $j$ the
identity component $G_j=\Aut_0(D_j)$ is noncompact. Write $G=\Aut_0(D)$.
Fix a $z_0\in D$. The orbit $G(z_0)$ is compact. We may
assume that $G(z_0)\sbs D_j$ for each $j$. For each connected
component $H$ of $\Aut(D)$, either the set $H(z_0)$ coincides with
$G(z_0)$ or $G(z_0)\cap H(z_0)=\emptyset$. Indeed, if $h\in H$ and
$h(z_0)\in G(z_0)$, then $H(z_0)=Gh(z_0)=G(z_0)$, since $H=Gh$. Now
we claim that there exists a positive number $a$ such that $a<
d(H(z_0), G(z_0))$ for each component $H$ of $\Aut(D)$ with
$H(z_0)\ne G(z_0)$, where $d$ is the euclidean distance. Otherwise, there is a sequence $\{H_k\}$ of
distinct components of $\Aut(D)$ with $H_k(z_0)\ne G(z_0)$ such that
$d(H_k(z_0), G(z_0))\to 0$. Passing to a subsequence if necessary, we
may assume that there are $h_k\in H_k$ such that  $h_k(z_0)$ tends to
a point in $G(z_0)$. It follows that some subsequence of $\{h_k\}$
converges in the compact-open topology to a $g\in \Aut(D)$; but this
is impossible because $h_k$ belong to different components of the Lie
group $\Aut(D)$. Therefore, such an $a$ exists. Decreasing $a$ if
necessary, we see that the open set
$$V=\{z\in D: d(z, G(z_0))<a\}$$
is relatively compact in $D$ and in each $D_j$, and satisfies $\overline
V\cap \Aut(D)(z_0)=G(z_0)$. This implies that $\bd V\cap
\Aut(D)(z_0)=\emptyset$. Since $G_j$ is noncompact, $G_j(z_0)$ is
noncompact, hence $G_j(z_0)\cap \bd V\ne \emptyset$. It follows that
for each $j$ there is a $g_j\in G_j$ with $g_j(z_0)\in \bd V$. Some
subsequence of the sequence $\{g_j\}$ converges uniformly on compacta
to a $g\in \Aut(D)$. It is clear that $g(z_0)\in \bd V$,
contradicting $\bd V\cap \Aut(D)(z_0)=\emptyset$.
\end{pf}


\begin{thebibliography}{999}
\bibitem[BD]{BD} E. Bedford, J. Dadok, {\em Bounded domains with
prescribed group of automorphisms,} Comment. Math. Helv. {\bf 62} (1987), 561--572
\bibitem[GK]{GK} R. Greene, S.G. Krantz, {\em Normal families and the
semicontinuity of isometry and automorphism groups,} Math. Z. {\bf
190} (1985), 455--467
\bibitem[FP]{FP} B. L. Fridman, E. A. Poletsky, {\sl Upper
semicontinuity of automorphism groups}, Math. Ann., 299(1994),
615--628
\bibitem[Ko]{Ko} Sh. Kobayashi, {\sl Hyperbolic Complex Spaces},
Springer-Verlag, 1998.
\bibitem[Ma]{Ma} D. Ma, {\sl Upper semicontinuity of isotropy and
automorphism groups}, Math.\ Ann., 292(1992), 533-545
\bibitem[MZ]{MZ} D. Montgomery, L. Zippin, {\em Topological transformation
groups}, New York, Interscience, 1955.
\bibitem[Sh]{Sh} B. V. Shabat,  {\em Introduction to Complex Analysis,
II}, AMS, 1992.
\bibitem[ShT]{ShT} A. E. Tumanov, G. B. Shabat, {\em Realization of
linear Lie groups by biholomorphic automorphisms of bounded domains,}
Funct. Anal. Appl. {\bf 24} (1990), 255--257
\bibitem[SZ]{SZ} R. Saerens, W. R. Zame, {\em The isometry groups of
manifolds and the automorphism groups of domains,} Trans. Amer. Math.
Soc. {\bf 301} (1987), 413--429
\end{thebibliography}
\end{document}